\documentclass[11pt, oneside, a4paper]{article}
\usepackage{amssymb}
\usepackage{color}
\usepackage{ulem}
\newcommand{\vp}{{\mathbf p}}
\newcommand{\vx}{{\mathbf x}}
\newcommand{\vz}{{\mathbf z}}
\newcommand{\vy}{{\mathbf y}}
\newcommand{\vu}{{\mathbf u}}
\newcommand{\vw}{{\mathbf w}}
\newcommand{\vf}{{\mathbf f}}

\newcommand{\la}{{\lambda}}

\newcounter{theorem}

\newcounter{theoremcounter}

\newcounter{remarkcounter}

\newcounter{definitioncounter}

\newtheorem{theorem}[theoremcounter]{Theorem}

\newtheorem{remark}[remarkcounter]{Remark}

\newtheorem{definition}[definitioncounter]{Definition}

\begin{document}



\title{On Obtaining Sharp Bounds  of the Rate of Convergence  for a Class of Continuous-Time Markov Chains}


\author{A. I. Zeifman\footnote{Vologda State University;  Institute of Informatics Problems,
Federal Research Center ``Computer Science and Control'',  Russian
Academy of Sciences; Vologda Research Center of the Russian Academy
of Sciences; e-mail a$\_$zeifman@mail.ru}; Y. A.
Satin\footnote{Vologda State University; e-mail yacovi@mail.ru}; K.
M. Kiseleva\footnote{Vologda State University; e-mail
ksushakiseleva@mail.ru}}

\date{}

\maketitle

{\bf Abstract.} We study inhomogeneous continuous-time weakly
ergodic Markov chains
 with a finite state space. We introduce
the notion of a Markov chain with the regular structure of
an infinitesimal matrix and study the sharp upper bounds on the rate of convergence for
such class of Markov chains.



\section{Introduction}
As is was noted in (Zeifman and Korolev, 2015), the problem of finding sharp
bounds on the rate of convergence to the limiting characteristics
for Markov chains is important for at least two following reasons:

(i) usually, it is easier to calculate the  limiting characteristics of
a process than to find the exact distribution of  its state probabilities.
Therefore, it is very important to have a possibility
to replace the exact distribution by the limiting one (see, for
instance, (Chakravarthy, 2017; Di Crescenzo et al. 2016)).
 In order to be able to determine the time
 instant starting from which such  replacement is possible, one must have
sharp bounds for the rate of convergence;

(ii) perturbation bounds play significant role in applications. It
is well known that even for homogeneous chains the best bounds
require the corresponding best bounds on the rate of convergence
(Kartashov, 1985, 1996;  Liu, 2012; Mitrophanov, 2003, 2004). Bounds
for general inhomogeneous Markov chains are also based on the
estimates of the rate of convergence in the special weighted norms
(Zeifman and Korolev, 2014).

In (Zeifman and Korolev, 2015) and (Zeifman et al. 2018b) we have studied
upper and lower bounds for finite inhomogeneous continuous-time
Markov chains. Moreover, in (Zeifman and Korolev, 2015) the new
approach to finding sharp upper bounds in natural metrics
via essential positivity of the reduced
intensity matrix of a Markov chain has been  presented.
 In this paper we
determine the class of  Markov chains for which this condition (i.e. essential positivity of the reduced intensity matrix; see matrix $B^*(t)$ below) is satisfied.
 It is shown  that this class is  fairly wide.

\medskip

Let $X(t)$ be  a finite continuous-time Markov chain with
infinitesimal matrix $Q(t)=\left(q_{i,j}(t)\right)_{i,j=0}^S$ and $A(t)=Q^T(t)$ - the transposed infinitesimal matrix of the process, respectively all $a_{i,j}(t)=q_{j,i}(t)$. Denote by $p_{ij}(s,t)=P\left\{ X(t)=j\left| X(s)=i\right.
\right\}$, $i,j \ge 0, \;0\leq s\leq t$ the transition probabilities
of $X(t)$ and by  $p_i(t)=P \left\{ X(t) =i \right\}$ -- the
probability  that the Markov chain $X(t)$ is in state $i$ at time
$t$. Let $\vp(t) = \left(p_0(t), p_1(t), \dots, p_S(t)\right)^T$ be the
probability distribution vector at instant $t$.

\begin{definition}
A Markov chain $X(t)$ is called a chain with the  regular infinitesimal
matrix if $q_{i,i+k}(t)$ and $q_{i+k,i}(t)$ decrease in $k$ for any
$t \ge 0$ and any fixed~$i$.
\end{definition}

\medskip

Consider the following four classes of finite inhomogeneous continuous-time Markov chains
 (see also (Zeifman et al. 2018a) and (Zeifman et al. 2018c) for more detailed description).

\medskip

\noindent (I) \textit{inhomogeneous birth-death processes}, where all $a_{ij}(t)=0$ for any $t\ge 0$ if $|i-j|>1$,
and $a_{i,i+1}(t)=\mu_{i+1}(t),~~a_{i+1,i}(t)=\la_i(t)$ - birth and death rates respectively;

\medskip

\noindent (II) \textit{inhomogeneous chains with ``batch'' births and single deaths},
where $a_{ij}(t)=0$ for any $t\ge 0$ if $i<j-1$, all ``birth'' rates do not depend on
the size of a ``population'', where $a_{i+k,i}(t)=a_k(t)$ for $k\ge 1$ - the rate of ``birth''
of a group of $k$ particles, $a_{i,i+1}(t)=\mu_{i+1}(t)$ - the death rate; under the
additional condition $a_{k+1}(t) \le a_k(t)$ for any positive $k$;

\medskip

\noindent (III) \textit{inhomogeneous chains with ``batch'' deaths and single births},
where all $a_{ij}(t)=0$ for any $t\ge 0$ if $i>j+1$, and all ``death'' rates do not
depend on the size of a ``population'', where $a_{i,i+k}(t)=b_k(t),~k\ge 1$ is the
rate of ``death'' of a group of $k$ particles, and $a_{i+1,i}(t)=\la_i(t)$ is the
birth rate; under the additional condition  $b_{k+1}(t) \le b_k(t)$ for any positive $k$;

\medskip

\noindent (IV) \textit{inhomogeneous chains with ``batch'' births and deaths}, where all
rates do not depend on the size of a ``population'', where $a_{i+k,i}(t)=a_k(t)$, and
$a_{i,i+k}(t)=b_k(t)$ for $k\ge 1$ are the rates of ``birth'' and ``death'' of a group
of $k$ particles respectively, under the additional condition
 $a_{k+1}(t) \le a_k(t)$ and $b_{k+1}(t) \le b_k(t)$, for any positive $k$.

\bigskip

One can see that {\it all} of the four Markov chains have the
regular structure of their infinitesimal matrices.

\section{Transformations for a class of Markov chains with the regular structure of an infinitesimal matrix}

Consider the forward Kolmogorov system  of ordinary differential equations for the
inhomogeneous Markov chain $X(t)$ in the form
\begin{equation} \label{ur01}
\frac{d}{dt} \vp(t)=A(t)\vp(t).
\end{equation}

 Using the substitution
$p_0(t) = 1 - \sum_{i= 1}^S p_i(t)$, we obtain from (\ref{ur01})
the reduced system in the form
\begin{equation}
\frac{d}{dt} \vz(t)= B(t)\vz(t)+\vf(t), \label{2.06}
\end{equation}
\noindent where $\vf(t)=\left(a_{10}(t), a_{20}(t),\cdots, a_{S0}(t)
\right)^T$, $\vz(t)=\left(p_{1}(t), p_{2}(t),\cdots, p_{S}(t)
\right)^T$,
\begin{equation}
{B(t) = \left(
\begin{array}{ccccccccc}
a_{11}(t)- a_{10}(t)   & a_{12}(t) - a_{10}(t)   &  \cdots & a_{1S}(t) - a_{10}(t)  \\
a_{21}(t) - a_{20}(t) & a_{22}(t) - a_{20}(t)   &    \cdots & a_{2S}(t) - a_{20}(t)  \\
a_{31}(t) - a_{30}(t)    & a_{32}(t) - a_{30}(t)  &    \cdots & a_{3S}(t) - a_{30}(t)  \\
\cdots \\
a_{S1}(t) - a_{S0}(t)  & a_{S2}(t) - a_{S0}(t) & \cdots     &
a_{SS}(t) -a_{S0}(t)
\end{array}
\right).} \label{2.07}
\end{equation}

All bounds on the rate of convergence to the limiting regime for
$X(t)$ correspond to the same   rate of convergence bounds of the solutions of system
\begin{equation}
\frac{d}{dt}\vy(t)= B(t)\vy(t). \label{hom1}
\end{equation}

Denote by $T$ the upper triangular matrix of the form
\begin{equation}
T=\left(
\begin{array}{ccccccc}
1   & 1 & 1 & \cdots & 1 \\
0   & 1  & 1  &   \cdots & 1 \\
0   & 0  & 1  &   \cdots & 1 \\
\vdots & \vdots & \vdots & \ddots \\
0   & 0  & 0  &   \cdots & 1
\end{array}
\right). \label{vspmatr}
\end{equation}

\bigskip

 Multiply (from the left) the left and the right part of (\ref{hom1}) by $T$ and
 put  $\vu(t)=T\vy(t)$.  In such a way the system (\ref{hom1})
is transformed into the system
\begin{equation}
\frac{d}{dt}\vu(t)= B^*(t)\vu(t), \label{hom11}
\end{equation}
\noindent where $B^*(t)=T
B(t)T^{-1}=\left(b^*_{ij}(t)\right)_{i,j=1}^S$
 and the entries $b^*_{ij}(t)$ having the following form:

$$b^*_{11}(t)=-a_{01}-a_{10}-a_{20}-\cdots-a_{S0},~b^*_{12}(t)=a_{01}-a_{02},$$
$$b^*_{13}(t)=a_{02}-a_{03},~ b^*_{14}(t)=a_{03}-a_{04}, ~\cdots ~, b^*_{1S}(t)=a_{0S-1}-a_{0S},$$
\medskip
$$b^*_{21}(t)=a_{21}-a_{20}+a_{31}-a_{30}+\cdots+a_{S1}-a_{S0},$$
$$b^*_{22}(t)=-a_{02}-a_{12}-a_{21}-a_{31}-\cdots-a_{S1},$$
$$b^*_{23}(t)=-a_{13}+a_{12}-a_{03}+a_{02},$$
$$b^*_{24}(t)=-a_{14}+a_{13}-a_{04}+a_{03},~\cdots,$$
$$b^*_{2S}(t)=-a_{1S}+a_{1S-1}-a_{0S}+a_{0S-1},$$
\medskip
$$b^*_{31}(t)=a_{31}-a_{30}+a_{41}-a_{40}+\cdots+a_{S1}-a_{S0},$$
$$b^*_{32}(t)=a_{32}-a_{31}+a_{42}-a_{41}+\cdots+a_{S2}-a_{S1},$$
$$b^*_{33}(t)=-a_{03}-a_{13}-a_{23}-a_{32}-a_{42}-\cdots-a_{S2},$$
$$b^*_{34}(t)=a_{03}-a_{04}+a_{13}-a_{14}+a_{23}-a_{24},~\cdots,$$
$$b^*_{3S}(t)=a_{0S-1}-a_{0S}+a_{1S-1}-a_{1S}+a_{2S-1}-a_{2S},$$
\medskip
$$b^*_{41}(t)=a_{41}-a_{40}+a_{51}-a_{50}+\cdots+a_{S1}-a_{S0},$$
$$b^*_{42}(t)=a_{42}-a_{41}+a_{52}-a_{51}+\cdots+a_{S2}-a_{S1},$$
$$b^*_{43}(t)=a_{43}-a_{42}+a_{53}-a_{52}+\cdots+a_{S3}-a_{S2},$$
$$b^*_{44}(t)=-a_{04}-a_{14}-a_{24}-a_{34}-a_{43}-a_{53}-\cdots-a_{S3},~\cdots,$$
$$b^*_{4S}(t)=a_{0S-1}-a_{0S}+a_{1S-1}-a_{1S}+a_{2S-1}-a_{2S}+a_{3S-1}-a_{3S},$$
\bigskip
$$b^*_{S1}(t)=a_{S1}-a_{S0},$$
$$b^*_{S2}(t)=a_{S2}-a_{S1},$$
$$b^*_{S3}(t)=a_{S3}-a_{S2},$$
$$b^*_{S4}(t)=a_{S4}-a_{S3}, \cdots, $$
$$b^*_{S,S}(t)=-a_{S,S-1}-a_{0S}-a_{1S}-a_{2S}-a_{3S}-\cdots-a_{S-1,S}.$$

\bigskip

 By inspecting the values of $b^*_{ij}(t)$
we arrive at the following theorem.
\begin{theorem} Let $X(t)$ be a finite Markov chain with the regular structure of
an infinitesimal matrix.  Then the corresponding respective matrix
$B^*(t)=T B(t)T^{-1}$ is essentially non-negative, i.e., $b^*_{ij} (t)
\ge0$ for all $i,j,$ such that $i\ne j$ and any $t \ge 0$.
\end{theorem}

\medskip

\begin{remark} One can note that the assumption of regularity of
infinitesimal matrix is close to the necessity for essential
non-negativity of $B^*(t)$. However, there are examples for which
infinitesimal matrix is not regular and $B^*(t)$ is essential
non-negative, see (Ammar and Alharbi, 2018; Dharmaraja, 2000; Zeifman et al. 2019).
\end{remark}

\section{General bounds}

In addition to the transformation $T$, described above,
we introduce  the new diagonal transformation
$D=diag\left(d_1,\dots,d_S\right)$, in which all the entries $d_j$ are positive.
 By left-multiplying the left and the right part of (\ref{hom11}) by $D$
and putting $\vw(t)=D\vu(t) = DT\vy(t)$, we obtain
\begin{equation}
\frac{d\vw}{dt}= B^{**}(t)\vw(t), \label{hom111}
\end{equation}
\noindent where $B^{**}(t)=D
B^*(t)D^{-1}=DT
B(t)T^{-1}D^{-1}=\left(b^{**}_{ij}(t)\right)_{i,j=1}^S$ is also essentially non-negative for any $t \ge 0$.

\medskip

Consider column sums of elements of $B^{**}(t)$ and the largest and smallest of them:
\begin{equation}
h^{**}(t)=\max_{1 \le k \le S}\sum_i b^{**}_{ik}(t), \quad h_{**}(t)=\min_{1 \le k \le S}\sum_i b^{**}_{ik}(t). \label{maxmin1}
\end{equation}

\medskip

 Using the same approach, which was used in the Theorem 1 of (Zeifman and Korolev, 2015),
we arrive at the following theorem.

\medskip

 \begin{theorem}  Let $X(t)$ be a finite Markov chain with the regular structure of
an infinitesimal matrix. Then for any positive $d_k$ the following bounds on the rate of convergence hold:
\begin{equation}
\|\vw(t)\| \le \exp\left\{{\int_0^t
h^{**}(\tau)\,d\tau}\right\}\|\vw(0)\|,
\label{301}
\end{equation}
\noindent for any corresponding initial condition $\vw(0)$, and
\begin{equation}
\|\vw(t)\| \ge \exp\left\{\int_0^t
h_{**}(\tau)\,d\tau\right\}\|\vw(0)\|,
\label{302}
\end{equation}
\noindent if the initial condition is non-negative,
\smallskip
where $\|\cdot\|$  denotes  the $l_1$-norm, i.e.,  $\|{\vx}\|=\sum|x_i|$.
\end{theorem}

\medskip

 In the homogeneous case,
sharp upper bounds on the rate of convergence can be obtained for
all of the fours classes of the considered Markov chains (classes (I)--(IV)).

\medskip

Firstly note that if a matrix $C$ is essentially non-negative and irreducible then
there exists a unique positive diagonal matrix $D= diag (d_1,
\dots, d_S)$ such that all column sums of matrix $C_D=DCD^{-1}$ are the
same (and are equal to the maximal eigenvalue $\lambda_0$ of $C$).
For the proof one can put  $m= \max{|c_{jj}|}$ and consider
non-negative matrix $C'=(C)^{T}+mI$. Matrix $C'$ is non-negative and
irreducible, hence it has a positive eigenvalue
$\lambda^*=\lambda_0+m$, and the respective eigenvector ${\bf x} =
\left(x_{1}, \dots, x_{S}\right)^T$ has positive coordinates. Put
$d_{k}=x_{k}^{-1}$. Then vector ${\bf e} = \left(1, \dots,
1\right)^T$ is eigenvector of matrix $C'_{D}=DC'D^{-1}$, and
therefore all row sums for this matrix is equal to $\lambda^*$.
Then all row sums for matrix $C^{T}_{D}= C'_{D} -mI$ is equal
to $\lambda^* -m=\lambda_0$, and every column sum for matrix $C_{D}$
will be the same. This implies an existence of positive diagonal
matrix $D= diag (d_1, \dots, d_S)$, and its uniqueness is evident.

\bigskip

\begin{theorem}  Let $X(t)$ be a finite homogeneous Markov chain with the regular structure of
an infinitesimal matrix. Then there exists a positive sequence
$\{d_k\}$ such that $h^{**}=h_{**}$, and hence
\begin{equation}
\|\vw(t)\| = \exp\left\{ h_{**}t\right\}\|\vw(0)\|, \label{3020}
\end{equation}
\noindent if the initial condition is non-negative,  and if:

\noindent (I) $X(t)$ is a birth-death processes with positive rates
$\mu_{i+1}$ and $\la_i$;

\noindent (II) $X(t)$ is a Markov chain with 'batch' births and
single deaths with positive $\mu_{i+1}$, and $a_{2} < a_1$;

\noindent (III) $X(t)$ is a Markov chain with   'batch' deaths and
single births with positive $\la_i$ and $b_{2} < b_1$;

\noindent (IV) $X(t)$ is a Markov chain with 'batch' births and
deaths with $a_{2} < a_1$ and $b_{2} < b_1$.

\end{theorem}

{\bf Proof.} For the first class we have
$$B^* = {\small \left(
\begin{array}{ccccc}
-\left(\lambda_0+\mu_1\right)  & \mu_1
 & 0 & \cdots & 0 \\
\lambda_1  & -\left(\lambda_1+\mu_2\right) & \mu_2 & \cdots & 0 \\
\ddots & \ddots & \ddots & \ddots & \ddots  \\
0 & \cdots & \cdots & \lambda_{S-1} &
-\left(\lambda_{S-1}+\mu_S\right)
\end{array}
\right).}$$

For the second class we obtain
$$B^* = {\small \left(
\begin{array}{ccccc}
a_{11}-a_S  & \mu_1
 & 0 & \cdots & 0 \\
a_1-a_S  & a_{22}-a_{S-1} & \mu_2 & \cdots & 0 \\
\ddots & \ddots & \ddots & \ddots & \ddots  \\
a_{S-1}-a_S & \cdots & \cdots & a_1-a_2 & a_{SS}-a_1
\end{array}
\right)}.$$

For the third class we have
$$B^* = {\small \left(
\begin{array}{ccccc}
-\left(\lambda_0 +b_1\right)  & b_1 - b_2
 & b_2 - b_3 & \cdots & b_{S-1} - b_S \\
\lambda_1  & -\big(\lambda_1+\sum\limits_{i\le 2}b_i\big) & b_1 - b_3 & \cdots & b_{S-2} - b_S \\
\ddots & \ddots & \ddots & \ddots & \ddots  \\
0 & \cdots & \cdots & \lambda_{S-1} &
-\big(\lambda_{S-1}+\sum\limits_{i\le S}b_i\big)
\end{array}
\right)}.$$

Finally,  for the fourth class we have
$$B^* = {\small \left(
\begin{array}{ccccc}
a_{11}-a_S  &  b_1 - b_2
 & b_2 - b_3 & \cdots & b_{S-1} - b_S\\
a_1-a_S  & a_{22}-a_{S-1} & b_1 - b_3 & \cdots & b_{S-2} - b_S \\ \\
\ddots & \ddots & \ddots & \ddots & \ddots  \\
a_{S-1}-a_S & \cdots & \cdots & a_1-a_2 & a_{SS}-a_1
\end{array}
\right)}.$$

\medskip

Note that in all cases,
all off-diagonal elements of the matrix $B^*$ are non-negative,
and all elements adjacent to the main diagonal are strictly positive.

\medskip

Hence under assumptions of  Theorem 3 all the matrices $B^*$ are essentially
non-negative and irreducible, hence all column sums for the corresponding
matrices $B^{**}$ are the same and $h^{**}=h_{**}$.

\medskip

\begin{remark} A finite birth-death process with the constant and
identical rates of birth $\lambda_k = a$ and $\mu_{k+1}=b$ was
considered in (Granovsky and Zeifman, 1997, 2005), where both upper
an lower sharp bonds were obtained. Namely, the corresponding
$D=D^*$, sharp $\beta_*=\beta^*= a+b-2\sqrt{ab}\cos\frac{\pi}{S+1}
\to (\sqrt{a}-\sqrt{b})^2$ as $S \to \infty$, and sharp $g^*=
a+b+2\sqrt{ab}\cos\frac{\pi}{S+1} \to (\sqrt{a}+\sqrt{b})^2$ as $S
\to \infty$ were found, where $\beta_*=\beta^* = -h^*=-h_*$.

\end{remark}

\medskip

\begin{remark} The possibility of obtaining explicit and exact estimates of the rate of convergence
 was previously studied for finite homogeneous birth-death processes, see
 (Granovsky and Zeifman, 1997, 2000; Van Doorn et al, 2010), as well as for chains of fourth class,
 see (Zeifman and Korolev, 2015).

 \end{remark}

\bigskip

{\bf Acknowledgement.} This research was supported by Russian
Science Foundation under grant 19-11-00020.

\newpage

\begin{center} {\bf References} \end{center}

\noindent Ammar, S. I.,  Alharbi, Y. F. 2018. Time-dependent
analysis for a two-processor heterogeneous system with time-varying
arrival and service rates. Applied Mathematical Modelling, 54, 743-751.

\medskip

\noindent Chakravarthy, S. R. 2017. A catastrophic queueing model
with delayed action. Applied Mathematical Modelling, 46, 631--649.

\medskip

\noindent Dharmaraja, S. 2000. Transient solution of a two-processor heterogeneous system, Math. Comput. Model. 32, 1117-–1123.

\medskip

\noindent Di Crescenzo, A., Giorno, V.,  Nobile, A. G. 2016.
Constructing transient birth-death processes by means of suitable
transformations. Applied Mathematics and Computation, 281, 152--171.

\medskip

\noindent Granovsky, B. L.,  Zeifman, A. I. 1997.  The decay function
of nonhomogeneous birth-death processes, with application to
mean-field models. Stochastic Process. Appl. 72, 105--120.

\medskip

\noindent Granovsky, B. L.,  Zeifman, A. I. 2000. The N-limit of spectral gap of a class
of birth–death Markov chains. Appl. Stoch. Models Bus. Ind. 16(4), 235--248.

\medskip

\noindent Granovsky, B. L.,  Zeifman, A. I. 2005. On the lower bound
of the spectrum of some mean-field models. Theory Probab. Appl. 49, 148--155.

\medskip

\noindent Kartashov, N.~V. 1985. Criteria for uniform ergodicity and
strong stability of Markov chains with a common phase space, Theory
Probab. Appl.  30, 71 -- 89.

\medskip

\noindent  Kartashov, N.V. 1996. Strong Stable Markov Chains. VSP. Kiev: TBiMC, Utrecht.

\medskip

\noindent  Liu, Y. 2012.  Perturbation bounds for the stationary distributions of Markov chains. SIAM. J. Matrix Anal. Appl. 33, 1057--1074.

\medskip

\noindent Mitrophanov, A.~Yu. 2003. Stability and exponential
convergence of continuous-time Markov chains. J. Appl. Probab.
40, 970--979.

\medskip

\noindent Mitrophanov, A.~Yu. 2004. The spectral gap and
perturbation bounds for reversible continuous-time Markov chains. J.
Appl. Probab.  41, 1219--1222.

\medskip

\noindent Van Doorn, E. A., Zeifman, A. I.,  Panfilova, T. L. 2010.
Bounds and asymptotics for the rate of convergence of birth-death processes.
Th. Probab. Appl. 54, 97--113.

\medskip

\noindent  Zeifman, A. I.,  Korolev, V. Y. 2014. On perturbation
bounds for continuous-time Markov chains. Statistics \& Probability
Letters, 88, 66--72.

\medskip

\noindent Zeifman, A. I.,  Korolev, V. Y. 2015. Two-sided Bounds on
the Rate of Convergence for Continuous-time Finite Inhomogeneous
Markov Chains. Statistics \& Probability Letters, 103, 30--36.

\medskip

\noindent Zeifman, A., Sipin, A., Korolev, V., Shilova, G., Kiseleva, K., Korotysheva, A.,  Satin, Y. 2018a. On  Sharp Bounds on the Rate of Convergence for Finite Continuous-time Markovian Queueing Models, LNCS, 10672, 20--28.

\medskip

\noindent Zeifman, A. I., Korolev, V. Y., Satin, Y. A.,  Kiseleva, K. M.
2018b. Lower bounds for the rate of convergence for continuous-time
inhomogeneous Markov chains with a finite state space. Statistics \&
Probability Letters, 137, 84--90.

\medskip

\noindent Zeifman, A., Razumchik, R., Satin, Y.,
Kiseleva, K., Korotysheva, A., Korolev, V. 2018c. Bounds on the rate of
convergence for one class of inhomogeneous Markovian queueing models
with possible batch arrivals and services. Int. J. Appl. Math. Comp.
Sci. 28, 141--154.

\medskip

\noindent Zeifman, A., Satin, Y., Kiseleva, K., Korolev, V.,  Panfilova, T. 2019. On limiting characteristics for a non-stationary two-processor heterogeneous system. Applied Mathematics and Computation, 351, 48--65.

\end{document}